\documentclass{article}

\usepackage{latexsym}

\usepackage{amsmath, amssymb}

\title{ZETA FUNCTIONS OF PERIODIC GRAPHS DERIVED FROM QUANTUM WALK}
\author{Takashi KOMATSU \\
Department of Bioengineering School of Engineering,\\ 
The University of Tokyo \\
Bunkyo, Tokyo, 113-8656, JAPAN \\ 
e-mail: komatsu@coi.t.u-tokyo.ac.jp \\ 
Norio KONNO \\
Department of Applied Mathematics, Faculty of Engineering, \\ 
Yokohama National University \\
Hodogaya, Yokohama 240-8501, JAPAN \\
e-mail: konno-norio-bt@ynu.ac.jp \\ 
Iwao SATO \\ 
Oyama National College of Technology \\
Oyama, Tochigi 323-0806, JAPAN \\ 
e-mail: isato@oyama-ct.ac.jp }
 
\begin{document}
 \maketitle

\clearpage

\vspace{5mm}



\vspace{5mm}

The contact author for correspondence: 

Iwao Sato 

Oyama National College of Technology, 
Oyama, Tochigi 323-0806, JAPAN

Tel: +81-285-20-2176

Fax: +81-285-20-2880

E-mail: isato@oyama-ct.ac.jp

\clearpage

\begin{abstract}
We define a zeta function of a finite graph derived from time evolution matrix of quantum walk, 
and give its determinant expression.  
Furthermore, we generalize the above result to a periodic graph. 
\end{abstract}

\section{Introduction}

Starting from $p$-adic Selberg zeta functions, Ihara [10] introduced the Ihara zeta functions 
of graphs. 
Ihara [10] showed that the reciprocal of the Ihara zeta function of a regular graph 
is an explicit polynomial. 
Serre [18] pointed out that the Ihara zeta function is 
the zeta function of a regular graph. 
A zeta function of a regular graph $G$ associated to a unitary 
representation of the fundamental group of $G$ was developed by 
Sunada [19,20]. 
Hashimoto [8] treated multivariable zeta functions of bipartite graphs. 
Bass [1] generalized Ihara's result on the Ihara zeta function of 
a regular graph to an irregular graph, and showed that its reciprocal is 
a polynomial. 

The Ihara zeta function of a finite graph was extended to an infinite graph 
in [1,3,5,6,7], and its determinant expressions were presented. 
Bass [1] defined the zeta function for a pair of a tree $X$ and a countable 
group $\Gamma $ which acts discretely on $X$ with quotient being a graph of 
finite groups. 
Clair and Mokhtari-Sharghi [3] extended Ihara zeta functions to infinite 
graphs on which a group $\Gamma $ acts isomorphically and with finite 
quotient. 
In [5], Grigorchuk and \.{Z}uk defined zeta functions of infinite discrete 
groups, and of some class of infinite periodic graphs. 
Guido, Isola and Lapidus [6] defined the Ihara zeta function of a 
periodic simple graph. 
Furthermore, Guido, Isola and Lapidus [7] presented a determinant expression for 
the Ihara zeta function of a periodic graph.

The time evolution matrix of a discrete-time quantum walk in a graph 
is closely related to the Ihara zeta function of a graph. 
A discrete-time quantum walk is a quantum analog of the classical random walk on a graph whose state vector is governed by 
a matrix called the time evolution matrix [12,16,21]. 
Ren et al. [17] gave a relationship between the discrete-time quantum walk and the Ihara zeta function of a graph.  
Konno and Sato [13] obtained a formula of the characteristic polynomial of the Grover matrix 
by using the determinant expression for the second weighted zeta function of a graph. 
Recently, Komatsu, Konno and Sato [11] presented a determinant expression for the zeta function with respect to 
the time evolution matrix of a general coined quantum walk on a periodic graph.

In this paper, we define a zeta function of a periodic graph by using the time evolution matrix 
of a generalization of a general coined quantum walk on it, and present its determinant expression. 

In Section 2, we state a review for the Ihara zeta function of a finite graph and 
infinite graphs, i.e., a periodic simple graph, a periodic graph. 
In Section 3, we state about the Grover walk on a graph as a discrete-time quantum walk on a graph. 
In Section 4, we define a zeta function of a finite graph $G$ by using the time evolution matrix 
of a generalization of a general coined quantum walk on $G$, and present its determinant expression. 
Furthermore, we give the explicit formula for its characteristic polynomial of $G$, and so present its spectrum. 
In Section 5, we state the definition of a periodic graph. 
In Section 6, we review a determinant for bounded operators acting on an 
infinite dimensional Hilbert space and belonging to a von Neumann algebra 
with a finite trace. 
In Section 7, we present a determinant expression for the above zeta function of a periodic graph.

\section{The Ihara zeta function of a graph}

All graphs in this paper are assumed to be simple. 
Let $G$ be a connected graph with vertex set $V(G)$ and edge set $E(G)$, 
and let $R(G)= \{ (u,v),(v,u) \mid uv \in E(G) \} $ be the set of 
oriented edges (or arcs) $(u,v),(v,u)$ directed oppositely for 
each edge $uv$ of $G$. 
For $e=(u,v) \in R(G)$, $u=o(e)$ and $v=t(e)$ are called 
the {\em origin} and the {\em terminal} of $e$, respectively. 
Furthermore, let $e^{-1}=(v,u)$ be the {\em inverse} of $e=(u,v)$. 

A {\em path $P$ of length $n$} in $G$ is a sequence 
$P=(e_1, \cdots ,e_n )$ of $n$ arcs such that $e_i \in R(G)$,
$t( e_i )=o( e_{i+1} )(1 \leq i \leq n-1)$. 
If $e_i =( v_{i-1} , v_i ), \  1 \leq i \leq n$, then we also denote 
$P$ by $( v_0, v_1 , \cdots ,v_n )$. 
Set $ \mid P \mid =n$, $o(P)=o( e_1 )$ and $t(P)=t( e_n )$. 
Also, $P$ is called an {\em $(o(P),t(P))$-path}. 
A $(v, w)$-path is called a {\em $v$-closed path} if $v=w$. 
The {\em inverse} of a closed path $C=( e_1, \cdots ,e_n )$ is the closed 
path $C^{-1} =( e^{-1}_n , \cdots ,e^{-1}_1 )$. 

We say that a path $P=(e_1, \cdots ,e_n )$ has a {\em backtracking} 
if $ e^{-1}_{i+1} =e_i $ for some $i(1 \leq i \leq n-1)$. 
A path without backtracking is called {\em proper}. 
Let $B^r$ be the closed path obtained by going $r$ times around a closed 
path $B$. 
Such a closed path is called a {\em multiple} of $B$. 
Multiples of a closed path without backtracking may have a backtracking. 
Such a closed path is said to have a {\em tail}. 
If its length is $n$, then the closed path can be written as 
\[
( e_1 , \cdots , e_k , f_1 , f_2 , \cdots , f_{n-2k} , 
e^{-1}_k , \cdots , e^{-1}_1 ) ,
\]
where $( f_1 , f_2 , \cdots , f_{n-2k} )$ is a closed path. 
A closed path is called {\em reduced} if $C$ has no backtracking nor tail. 
Furthermore, a closed path $C$ is {\em primitive} if it is not a multiple of 
a strictly shorter closed path. 

We introduce an equivalence relation between closed paths. 
Two closed paths $C_1 =(e_1, \cdots ,e_m )$ and 
$C_2 =(f_1, \cdots ,f_m )$ are called {\em equivalent} if 
there exists an integer $k$ such that $f_j =e_{j+k} $ for all $j$, 
where the subscripts are read modulo $n$. 
The inverse of $C$ is not equivalent to $C$ if $\mid C \mid \geq 3$. 
Let $[C]$ be the equivalence class which contains a closed path $C$. 
Also, $[C]$ is called a {\em cycle}. 

Let ${\cal P} $ be the set of primitive, reduced cycles of $G$.  
Also, primitive, reduced cycles are called {\em prime cycles}. 
Note that each equivalence class of primitive, reduced closed paths 
of a graph $G$ passing through a vertex $v$ of $G$ corresponds 
to a unique conjugacy class of the fundamental group $ \pi {}_1 (G,v)$ 
of $G$ at $v$. 

The {\em Ihara zeta function} of a graph $G$ is 
a function of a complex variable $u$ with $|u|$ 
sufficiently small, defined by 
\[
{\bf Z} (G, u)= {\bf Z}_G (u)= \prod_{[C] \in {\cal P}} 
(1- u^{ \mid C \mid } )^{-1} ,
\]
where $[C]$ runs over all prime cycles of $G$.

Let $G$ be a connected graph with $n$ vertices $v_1, \cdots ,v_n $. 
The {\em adjacency matrix} ${\bf A}= {\bf A} (G)=(a_{ij} )$ is 
the square matrix such that $a_{ij} =1$ if $v_i$ and $v_j$ are adjacent, 
and $a_{ij} =0$ otherwise.
The {\em degree} of a vertex $v_i$ of $G$ is defined by 
$ \deg v_i = \deg {}_G v_i = \mid \{ v_j \mid v_i v_j \in E(G) \} \mid $. 
If $ \deg {}_G v=k$(constant) for each $v \in V(G)$, then $G$ is called 
{\em $k$-regular}.

\newtheorem{theorem}{Theorem}
\begin{theorem}[Bass] 
Let $G$ be a connected graph. 
Then the reciprocal of the Ihara zeta function of $G$ is given by 
\[
{\bf Z} (G,u )^{-1} =(1- u^2 )^{r-1} 
\det ( {\bf I} -u {\bf A} (G)+ u^2 ( {\bf D} - {\bf I} )) , 
\]
where $r$ is the Betti number of $G$, 
and ${\bf D} =( d_{ij} )$ is the diagonal matrix 
with $d_{ii} = \deg v_i$ and $d_{ij} =0, i \neq j , 
(V(G)= \{ v_1 , \cdots , v_n \} )$. 
\end{theorem}

Let $G=(V(G),E(G))$ be a countable simple graph, 
and let $\Gamma $ be a countable discrete subgroup of 
automorphisms of $G$, which acts freely on $G$, and with finite quotient 
$G/ \Gamma $. 
The graph $G$ is called a {\em periodic graph}.  
Then the Ihara zeta function of a periodic simple graph is defined 
as follows: 
\[
{\bf Z}_{G, \Gamma} (u)= \prod_{[C]_{\Gamma} \in [{\cal P}]_{\Gamma }} 
(1- u^{\mid C \mid } )^{- 1/ \mid \Gamma {}_{[C]} \mid } , 
\]
where $ \Gamma {}_{[C]} $ is the stabilizer of $[C]$ in $\Gamma $, and 
$[C]_{\Gamma} $ runs over all $\Gamma $-equivalence classes of 
prime cycles in $G$. 

Guido, Isola and Lapidus [6] presented a determinant expression for 
the Ihara zeta function of a periodic simple graph.

\begin{theorem}[Guido, Isola and Lapidus] 
For a periodic simple graph $G$, 
\[
{\bf Z}_{G, \Gamma } (u)=(1- u^2 )^{-(m-n)} 
\det {}_{\Gamma } ( {\bf I} -u {\bf A} (G) +( {\bf D} - {\bf I} ) 
u^2 )^{-1} , 
\]
where $\det {}_{\Gamma } $ is a determinant for bounded operators 
belonging to a von Neumann algebra with a finite trace. 
\end{theorem}

Guido, Isola and Lapidus [7] presented a determinant expression for 
the Ihara zeta function of a periodic graph $G$ and a countable discrete subgroup 
$\Gamma $ of aoutomorphisms of $G$ which acts discretely without inversions, and 
with bounded covolume.

\begin{theorem}[Guido, Isola and Lapidus] 
For a periodic graph $G$, 
\[
{\bf Z}_{G, \Gamma } (u)^{-1} =(1- u^2 )^{\chi {}^{(2)} (G)} 
\det {}_{\Gamma } ( \Delta (u)) , 
\]
where ${\chi {}^{(2)} (G)}$ is the $L^2$-Euler characteristic of $(G, \Gamma )$
(see [2]), and $\Delta (u)={\bf I} -u{\bf A}+ u^2 ( {\bf D} -{\bf I} )$. 
\end{theorem}

\section{The general coined quantum walk on a graph}

First, we state the definition of a coined quantum walk as a definition of 
a discrete-time quantum walk on a graph(see [12,16,17]). 

Let $G$ be a connected graph with $m$ edges. 
For each arc $e=(u,v) \in D(G)$, we indicate the {\em pure state} $|e \rangle =|uv \rangle $ such that 
$\{ |e \rangle \mid e \in D(G) \} $ is a normal orthogonal system on the Hilbert space $\mathbb{C}^{2m} $. 
The {\em transition} from an arc $(u,v)$ to an arc $(w,x)$ occurs if $v=w$. 
The {\em state} of quantum walk is defined as follows: 
\[
\psi = \sum_{(u,v) \in D(G)} \alpha {}_{uv} |uv \rangle , \  \alpha {}_{uv} \in \mathbb{C} . 
\]
The {\em probability} that there exists a particle in the arc $e=(u,v)$ is given as follows: 
\[
P( |e \rangle )= \alpha {}_{uv} \overline{ \alpha {}_{uv}} . 
\]
Here, 
\[
\sum_{(u,v) \in D(G)} \alpha {}_{uv} \overline{ \alpha {}_{uv}} =1 . 
\]

Let $G$ be a connected graph with $n$ vertices and $m$ edges. 
Set $V(G)= \{ v_1 , \ldots , v_n \} $ and $d_j = d_{v_j} = \deg v_j , \ j=1, \ldots , n$. 
For $u \in V(G)$, let $D(u)= \{ e \in D(G) \mid o(e)=u \} $. 
Then, for $u \in V(G)$, let 
\[
D(u)= \{ e_{u1} , \ldots , e_{u d_u } \} .
\] 
Furthermore, let $\alpha {}_u , \ u \in V(G)$ be a unit vector with respect to $D(u)$, 
that is, 
\[
\alpha {}_u (e) =\left\{
\begin{array}{ll}
non \ zero \ complex \ number \ & \mbox{if $e \in D(u)$, } \\
0 & \mbox{otherwise, }
\end{array}
\right.
\]
where $ \alpha {}_u (e) $ is the entry of $\alpha {}_u $ corresponding to the arc $e \in D(G)$. 

Now, a $2m \times 2m$ matrix ${\bf C} $ is given as follows: 
\[
{\bf C} =2 \sum_{u \in V(G)} | \alpha {}_u \rangle \langle \alpha {}_u | - {\bf I}_{2m } . 
\]
The matrix {\bf C} is the {\em coin  operator} of the considered quantum walk. 
Note that ${\bf C} $ is unitary. 
Then the {\em transition matrix} ${\bf U} $ is defined by 
\[
{\bf U} = {\bf J}_0 {\bf C} . 
\]
The matrix ${\bf J}_0 $ is called the {\em shift operator}.

The time evolution of a quantum walk on $G$ through ${\bf U} $ is given by 
\[
\psi {}_{t+1} = {\bf U} \psi {}_t . 
\]
Here, $\psi {}_{t+1}, \psi {}_t $ are the states. 
Note that the state $\psi {}_t $ is written with respect to the initial state $\psi {}_0 $ as follows: 
\[ 
\psi {}_{t} = {\bf U}^t \psi {}_0 . 
\]
A quantum walk on $G$ with ${\bf U} $ as a time evolution matrix is called a 
{\em coined quantum walk} on $G$.

We consider a general coined quantum walk on a graph. 
We replace the coin operator ${\bf C} $ of a coined quantum walk with unitary matrix with two spectra which are distinct 
from $\pm 1$. 

For a given connected graph $G$ with $n$ vertices and $m$ edges, let ${\bf d} : \ell {}^2 (V(G)) \longrightarrow \ell {}^2 (R(G))$ 
such that 
\[
{\bf d} {\bf d}^* ={\bf I}_{q} ,  
\] 
and let ${\bf S} =( S_{ef} )_{e,f \in R(G)} $ be the $2m \times 2m$ matrix defined by 
\[
S_{ef} =\left\{
\begin{array}{ll}
1 & \mbox{if $f= e^{-1} $, } \\
0 & \mbox{otherwise.}
\end{array}
\right.
\] 
Then the matrices ${\bf S} $ and ${\bf C} $ are called the {\em shift operator} and the {\em coin operator} of ${\bf U} $, 
respectively.

Furthermore, let 
\[
{\bf C} = a {\bf d}^* {\bf d}+b( {\bf I}_{2m} - {\bf d}^* {\bf d} ) 
\] 
and ${\bf U} = {\bf S} {\bf C} $(see [9]). 
Note that $q= \dim \ ker(a- {\bf C} )$.
A discrete-time quantum walk on $G$ with ${\bf U}$ as a time evolution matrix is called a {\em general coined quantum walk} 
on $G$. 
Then we define a zeta function of $G$ by using ${\bf U}$ as follows: 
\[
\zeta {} (G, u)= \det ( {\bf I}_{2m} -u {\bf U} )^{-1} 
= \det ( {\bf I}_{2m} -u {\bf S} (a {\bf d}^* {\bf d} +b( {\bf I}_{2m} - {\bf d}^* {\bf d} )) )^{-1} . 
\]

Now, we have the following result.

\begin{theorem}[Komatsu, Konno and Sato]
Let $G$ be a connected graph $n$ vertices and $m$ edges, ${\bf U} = {\bf S} {\bf C} $ the time evolution matrix 
of a general coined quantum walk on $G$. 
Suppose that $\sigma ( {\bf C} )= \{ a,b \} $. 
Set $q= \dim \ ker(a- {\bf C} )$.
Then, for the unitary matrix ${\bf U} = {\bf S} {\bf C} $, we have    
\[
\zeta {} (G, u)=(1- b^2 u^2 )^{m-q} \det ((1-abu^2 ) {\bf I}_n -cu {\bf d} {\bf S} {\bf d}^* ) , c=a-b .   
\]
\end{theorem}

\section{Spectra for the unitary matrix of a generalization of a general coined quantum walk on a graph} 

We consider a generalization of a general coined quantum walk on a graph. 
We replace the coin operator ${\bf C} $ and the shift operator ${\bf S} $of the above quantum walk with unitary matrices with two spectra 
which are distinct from $\pm 1$. 

For given connected graph $G$ with $n$ vertices and $m$ edges, let 
${\bf d}_i : \ell {}^2 (V(G)) \longrightarrow \ell {}^2 (R(G)) \ (i=1,2)$ such that 
\[
{\bf d}_1 {\bf d}^*_1 ={\bf I}_{p}\   , {\bf d}_2 {\bf d}^*_2 ={\bf I}_{q} . 
\] 
Furthermore, let 
\[
{\bf C}_i = a_i {\bf d}^*_i {\bf d}_i +b_i ( {\bf I}_{2m} - {\bf d}^*_i {\bf d}_i ) \ (i=1,2)  
\] 
and ${\bf U} = {\bf C}_1 {\bf C}_2 $(see [9]). 
A discrete-time quantum walk on $G$ with ${\bf U}$ as a time evolution matrix is called a {\em general coined quantum walk} 
on $G$. 
Then we define a zeta function of $G$ by using ${\bf U}$ as follows: 
\[
\zeta {} (G, u)= \det ( {\bf I}_{2m} -u {\bf U} )^{-1} 
= \det ( {\bf I}_{2m} -u (a_1 {\bf d}^*_1 {\bf d}_1 +b_1 ( {\bf I}_{2m} - {\bf d}^*_1 {\bf d}_1 ))
(a_2 {\bf d}^*_2 {\bf d}_2 +b_2 ( {\bf I}_{2m} - {\bf d}^*_2 {\bf d}_2 )) )^{-1} . 
\]

Now, we have the following result.

\begin{theorem}
Let $G$ be a connected graph $n$ vertices and $m$ edges, ${\bf U} = {\bf C}_1 {\bf C}_2 $ the time evolution matrix of a coined quantum walk on $G$. 
Suppose that $\sigma ( {\bf C}_i )= \{ a_1 ,b_i \} \ (i=1,2)$. 
Set $p= \dim \ ker(a_1 - {\bf C}_1 )$ and $q= \dim \ ker(a_2 - {\bf C}_2 )$. 
Then, for the unitary matrix ${\bf U} = {\bf C\natural }_1 {\bf C}_2 $, we have    
\[
\zeta {} (G, u)=(1- b_1 b_2 u )^{2m-p-q} (1- a_2 b_1 u )^{q-p} \det ((1-a_1 b_2 u)(1-a_2 b_1 u) {\bf I}_q -c_1 c_2 u {\bf d}_1 {\bf d}^*_2 {\bf d}_2 {\bf d}^*_1 ) , 
\] 
where 
\[
c_i =a_i -b_i \ (i=1,2) .   
\]
\end{theorem}

{\em Proof }.  If ${\bf A}$ and ${\bf B}$ are a $m \times n $ and $n \times m$ 
matrices, respectively, then we have 
\[
\det ( {\bf I}_{2m} - {\bf A} {\bf B} )= 
\det ( {\bf I}_n - {\bf B} {\bf A} ) . 
\]    
Let 
\[
c_i = a_i -b_i \  (i=1,2) . 
\]
Then, for ${\bf U} = {\bf C}_1 {\bf C}_2 $, 
\[
\begin{array}{rcl} 
\  &   & \det ( {\bf I}_{2m} -u {\bf U} )= \det ( {\bf I}_{ 2m} -u {\bf C}_1 {\bf C}_2 ) \\ 
\  &   &                \\ 
\  & = & 
\det ({\bf I}_{2m} -u(c_1 {\bf d}^*_1 {\bf d}_1 +b_1 {\bf I}_{2m} )
( c_2 {\bf d}^*_2 {\bf d}_2 + b_2 {\bf I}_{2m } )) \\ 
\  &   &                \\ 
\  & = & \det ({\bf I}_{2m} -c_1 u {\bf d}^*_1 {\bf d}_1 (c_2 {\bf d}^*_2 {\bf d}_2 + b_2 {\bf I}_{2m } )
-b_1 u(c_2 {\bf d}^*_2 {\bf d}_2 + b_2 {\bf I}_{2m } )) \\ 
\  &   &                \\ 
\  & = & 
\det ((1- b_1 b_2 u) {\bf I}_{2m} -b_1 c_2 u {\bf d}^*_2 {\bf d}_2 -c_1 u {\bf d}^*_1 {\bf d}_1 (c_2 {\bf d}^*_2 {\bf d}_2 + b_2 {\bf I}_{2m } )) \\ 
\  &   &                \\ 
\  & = & (1-u )^{2m} 
\det ( {\bf I}_{2m} - \frac{b_1 c_2 u}{1-b_1 b_2 u} {\bf d}^*_2 {\bf d}_2  
- \frac{cu}{1-b_1 b2 u} {\bf d}^*_1 {\bf d}_1 (c_2 {\bf d}^*_2 {\bf d}_2 + b_2 {\bf I}_{2m } )) \\ 
\  &   &                \\ 
\  & = & (1-b_1 b_2 u )^{2m} 
\det ( {\bf I}_{2m } - \frac{c_1 u}{1-b_1 b_2 u} {\bf d}^*_1 {\bf d}_1 (c_2 {\bf d}^*_2 {\bf d}_2 + b_2 {\bf I}_{2m } )
( {\bf I}_{2m } - \frac{b_1 c_2 u}{1-b_1 b_2 u} {\bf d}^*_2 {\bf d}_2 )^{-1} ) \\ 
\  &   &                \\ 
\  & \times & \det ( {\bf I}_{2m } - \frac{b_1 c_2 u}{1-b_1 b_2 u} {\bf d}^*_2 {\bf d}_2 ) . 
\end{array}
\]

But, we have 
\[
\begin{array}{rcl}
\det ( {\bf I}_{2m } - \frac{b_1 c_2 u}{1-b_1 b_2 u} {\bf d}^*_2 {\bf d}_2 ) 
& = & \det ( {\bf I}_{q} - \frac{b_1 c_2 u}{1-b_1 b_2 u} {\bf d}_2 {\bf d}^*_2 ) \\
\  &   &                \\ 
\  & = & \det ( {\bf I}_{q} - \frac{b_1 c_2 u}{1-b_1 b_2 u} {\bf I}_q ) \\
\  &   &                \\ 
\  & = & (1- \frac{b_1 c_2 u}{1-b_1 b_2 u} )^q = \frac{(1-b_1 a_2 u)^q }{(1-b_1 b_2 u)^q } . 
\end{array}
\] 
Furthermore, we have 
\[
\begin{array}{rcl}
\  &   & ( {\bf I}_{2m } - \frac{b_1 c_2 u}{1-b_1 b_2 u} {\bf d}^*_2 {\bf d}_2 )^{-1} \\ 
\  &   &                \\ 
\  & = & {\bf I}_{2m } + \frac{b_1 c_2 u}{1-b_1 b_2 u} {\bf d}^*_2 {\bf d}_2 +(\frac{b_1 c_2 u}{1-b_1 b_2 u} )^2 {\bf d}^*_2 {\bf d}_2 {\bf d}^*_2 {\bf d}_2  
+ ( \frac{b_1 c_2 u}{1-b_1 b_2 u} )^3 {\bf d}^*_2 {\bf d}_2 {\bf d}^*_2 {\bf d}_2 {\bf d}^*_2 {\bf d}_2 + \cdots \\
\  &   &                \\ 
\  & = & {\bf I}_{2m } + \frac{b_1 c_2 u}{1-b_1 b_2 u} {\bf d}^*_2 {\bf d}_2 +( \frac{b_1 c_2 u}{1-b_1 b_2 u} )^2 {\bf d}^*_2 {\bf d}_2 
+ ( \frac{b_1 c_2 u}{1-b_1 b_2 u} )^3 {\bf d}^*_2 {\bf d}_2 + \cdots \\
\  &   &                \\ 
\  & = & {\bf I}_{2m } + \frac{b_1 c_2 u}{1-b_1 b_2 u} (1 + \frac{b_1 c_2 u}{1-b_1 b_2 u} +( \frac{b_1 c_2 u}{1-b_1 b_2 u} )^2 - \cdots ) {\bf d}^*_2 {\bf d}_2 \\
\  &   &                \\ 
\  & = & {\bf I}_{2m } + \frac{b_1 c_2 u}{1-b_1 b_2 u} /(1- \frac{b_1 c_2 u}{1-b_1 b_2 u} ) {\bf d}^*_2 {\bf d}_2 )
={\bf I}_{2m } + \frac{b_1 c_2 u}{1-b_1 a_2 u} {\bf d}^*_2 {\bf d}_2 .
\end{array}
\]

Therefore, it follows that 
\[
\begin{array}{rcl}
\  &   & \det ( {\bf I}_{2m } -u {\bf U} ) \\ 
\  &   &                \\ 
\  & = & (1-b_1 b_2 u)^{2m} 
\det ( {\bf I}_{2m } - \frac{c_1 u}{1-b_1 b_2 u} {\bf d}^*_1 {\bf d}_1 (c_2 {\bf d}^*_2 {\bf d}_2 + b_2 {\bf I}_{2m } )
( {\bf I}_{2m} + \frac{b_1 c_2 u}{1-b_1 a_2 u} {\bf d}^*_2 {\bf d}_2 )) \frac{(1-b_1 a_2 u)^q }{(1-b_1 b_2 u)^q } \\ 
\  &   &                \\ 
\  & = & (1-b_1 b_2 u)^{2m-q} (1-b_1 a_2 u)^q 
\det ( {\bf I}_{2m } - \frac{c_1 u}{1-b_1 b_2 u} {\bf d}^*_1 {\bf d}_1 (b_2 {\bf I}_{2m } + \frac{ c_2 }{1-b_1 a_1 u} {\bf d}^*_2 {\bf d}_2 )) \\ 
\  &   &                \\ 
\  & = & (1-b_1 b_2 u)^{2m-q} (1-b_1 a_2 u)^q  
\det ( {\bf I}_{p} - \frac{c_1 u}{1-b_1 b_2 u} {\bf d}_1 (b_2 {\bf I}_{2m } + \frac{ c_2 }{1-b_1 a_1 u} {\bf d}^*_2 {\bf d}_2 ) {\bf d}^*_1 ) \\ 
\  &   &                \\ 
\  & = & (1-b_1 b_2 u)^{2m-q} (1-b_1 a_2 u)^q     
\det ( {\bf I}_{p} - \frac{c_1 u}{1-b_1 b_2 u} {\bf d}_1 {\bf d}^*_1 - \frac{c_1 c_2 u}{(1-b_1 b_2 u )(1-b_1 a_2 u)} {\bf d}_1 {\bf d}^*_2 {\bf d}_2 {\bf d}^*_1 ) \\ 
\  &   &                \\ 
\  & = & (1-b_1 b_2 u)^{2m-q} (1-b_1 a_2 u)^q  
\det ( {\bf I}_{p} - \frac{c_1 u}{1-b_1 b_2 u} {\bf I}_p - \frac{c_1 c_2 u}{(1-b_1 b_2 u )(1-b_1 a_2 u)} {\bf d}_1 {\bf d}^*_2 {\bf d}_2 {\bf d}^*_1 ) \\ 
\  &   &                \\ 
\  & = & (1-b_1 b_2 u)^{2m-q} (1-b_1 a_2 u)^q   
\det ( \frac{1-a_1 b_2 u}{1-b_1 b_2 u} {\bf I}_{p} - \frac{c_1 c_2 u}{(1-b_1 b_2 u )(1-b_1 a_2 u)} {\bf d}_1 {\bf d}^*_2 {\bf d}_2 {\bf d}^*_1 ) \\ 
\  &   &                \\ 
\  & = & (1-b_1 b_2 u)^{2m-p-q} (1-b_1 a_2 u)^{q-p} \det ((1-a_1 b_2 u)(1-b_1 a_2 u) {\bf I}_{p} -c_1 c_2 u {\bf d}_1 {\bf d}^*_2 {\bf d}_2 {\bf d}^*_1) . 
\end{array}
\]
Q.E.D.

\newtheorem{corollary}{Corollary} 
\begin{corollary} 
Let $G$ be a connected graph $n$ vertices and $m$ edges, ${\bf U} = {\bf C}_1 {\bf C}_2 $ the time evolution matrix of a coined quantum walk on $G$. 
Suppose that $\sigma ( {\bf C}_i )= \{ a_1 ,b_i \} \ (i=1,2)$. 
Set $p= \dim \ ker(a_1 - {\bf C}_1 )$ and $q= \dim \ ker(a_2 - {\bf C}_2 )$.
Then, for the unitary matrix ${\bf U} = {\bf C\natural }_1 {\bf C}_2 $, we have    
\[
\det ( \lambda {\bf I}_{2m} - {\bf U} )=( \lambda - b_1 b_2 )^{2m-p-q} ( \lambda - a_2 b_1 u )^{q-p} 
\det (( \lambda -a_1 b_2 )( \lambda -a_2 b_1 ) {\bf I}_p -c_1 c_2 \lambda {\bf d}_1 {\bf d}^*_2 {\bf d}_2 {\bf d}^*_1 ) , 
\] 
where 
\[
c_i =a_i -b_i \ (i=1,2) .   
\]
\end{corollary}

{\em Proof }.  Let $u=1/ \lambda $. 
Then, by Theorem 6, we have 
\[
\det ({\bf I}_{2m} -1/  \lambda {\bf U} )=(1- b_1 b_2 / \lambda )^{2m-p-q} (1- a_2 b_1 / \lambda )^{q-p} 
\det ((1-a_1 b_2 / \lambda)(1-a_2 b_1 / \lambda) {\bf I}_p -c_1 c_2 / \lambda {\bf d}_1 {\bf d}^*_2 {\bf d}_2 {\bf d}^*_1 ) , 
\] 
and so, 
\[
\det ( \lambda  {\bf I}_{2m} - {\bf U} )=( \lambda - b_1 b_2 )^{2m-p-q} ( \lambda - a_2 b_1 u )^{q-p} 
\det (( \lambda -a_1 b_2 )( \lambda -a_2 b_1 ) {\bf I}_p -c_1 c_2 \lambda {\bf d}_1 {\bf d}^*_2 {\bf d}_2 {\bf d}^*_1 ) , 
\] 
$\Box$

By Corollary 1, the following result holds.

\begin{corollary} 
Let $G$ be a connected graph $n$ vertices and $m$ edges, ${\bf U} = {\bf C}_1 {\bf C}_2 $ the time evolution matrix of a coined quantum walk on $G$. 
Suppose that $\sigma ( {\bf C}_i )= \{ a_1 ,b_i \} \ (i=1,2)$. 
Set $p= \dim \ ker(a_1 - {\bf C}_1 )$ and $q= \dim \ ker(a_2 - {\bf C}_2 )$.
Then, the spectra of the unitary matrix ${\bf U} = {\bf C}_1 {\bf C}_2 $ are given as follows: 
\begin{enumerate} 
\item $2p$ eigenvalues: 
\[
\lambda = \frac{a_1 b_2 +b_1 a_2 +c_1 c_2 \mu \pm \sqrt{ ( a_1 b_2 +b_1 a_2 +c_1 c^2 \mu )^2 -4a_1 b_1 a_2 b_2 \mu }}{2} , 
\ \mu \in Spec ( {\bf d}_1 {\bf d}^*_2 {\bf d}_2 {\bf d}^*_1 ) ; 
\]
\item $q-p$ eigenvalues: $b_1 a_2 $ ; 
\item $2m-p-q$ eigenvalues: $b_1 b_2$ .  
\end{enumerate} 
\end{corollary}

{\em Proof }. By Corollary 1, we have 
\[
\begin{array}{rcl}
\  &   & \det ( \lambda  {\bf I}_{2m} - {\bf U} ) \\
\  &   &                \\ 
\  & = & ( \lambda - b_1 b_2 )^{2m-p-q} ( \lambda - a_2 b_1 u )^{q-p} 
\prod_{ \mu \in Spec ( {\bf d}_1 {\bf d}^*_2 {\bf d}_2 {\bf d}^*_1 )} (( \lambda -a_1 b_2 )( \lambda -a_2 b_1 )-c_1 c_2 \mu \lambda ) . 
\end{array}
\]

Solving 
\[
( \lambda -a_1 b_2 )( \lambda -a_2 b_1 )-c_1 c_2 \mu \lambda = \lambda {}^2 -(a_1 b_2 +b_1 a_2 +c_1 c^2 \mu ) \lambda +a_1 b_2 a_2 b_1 =o , 
\]
we obtain 
\[
\lambda = \frac{a_1 b_2 +b_1 a_2 +c_1 c_2 \mu \pm \sqrt{ ( a_1 b_2 +b_1 a_2 +c_1 c^2 \mu )^2 -4a_1 b_1 a_2 b_2 \mu }}{2} , 
\] 
The result follows. 
$\Box$

\section{Periodic graphs}

Let $G=(V(G),E(G))$ be a simple graph. 
Assume that $G$ is countable ($V(G)$ and $E(G)$ are countable), 
and with bounded degree, i.e., $d= \sup_{v \in V(G)} \deg v < \infty$. 
Let $\Gamma $ be a countable discrete subgroup of automorphisms of $G$, 
which acts 
\begin{enumerate} 
\item  without inversions: $\gamma (e) \neq e^{-1} $ for any 
$\gamma \in \Gamma, e \in R(G)$, 
\item  discretely: $\Gamma {}_v = \{ \gamma \in \Gamma \mid \gamma v=v \}$ 
is finite for any $v \in V(G)$, 
\item  with bounded covolume: ${\rm vol} (G/ \Gamma ):= 
\sum_{v \in {\cal F}_0 } \frac{1}{ \mid \Gamma {}_v \mid } < \infty $, where 
${\cal F}_0 \subset V(G)$ contains exactly one representative for 
each equivalence class in $V(G/ \Gamma )$. 
\end{enumerate} 
Then $G$ is called a {\em periodic graph} with a countable discrete subgroup 
$\Gamma $ of $Aut \  G$. 
Note that the third condition is equivalent to the following condition: 
\[
{\rm vol} (R(G)/ \Gamma ):= 
\sum_{e \in {\cal F}_1 } \frac{1}{ \mid \Gamma {}_e \mid } < \infty , 
\]
where ${\cal F}_1 \subset R(G)$ contains exactly one representative for 
each equivalence class in $R(G/ \Gamma )$. 

Let $\ell {}^2 (V(G))$ be the Hilbert space of functions 
$f: V(G) \longrightarrow {\bf C} $ such that 
$\mid \mid f \mid \mid := \sum_{v \in V(G)} \mid f(v) \mid {}^2 
< \infty $. 
We define the left regular representation $\lambda {}_0 $ of $\Gamma $ on 
$\ell {}^2 (V(G))$ as follows: 
\[
(\lambda {}_0 ( \gamma )f)(x)=f( \gamma {}^{-1} x), \  \gamma \in \Gamma, 
\  f \in \ell {}^2 (V(G)), \  x \in V(G) . 
\]

We state the definition of a von Neumann algebra. 
Let $H$ be a separable complex Hilbert space, and let ${\cal B} (H)$ denote 
the ${\bf C}^*$-algebra of bounded linear operators on $H$. 
For a subset $M \subset {\cal B} (H)$, the {\em commutant} of $M$ is 
$M^{\prime} = \{ T \in {\cal B} (H) \mid ST=TS, \forall S \in M \} $. 
Then a {\em von Neumann algebra} is a subalgebra ${\cal A} \leq {\cal B} (H)$ 
such that ${\cal A}^{\prime \prime} ={\cal A} $. 
It is known that a determinant is defined for a suitable class of operators 
in a von Neumann algebra with a finite trace (see [4,6]). 
 
For the Hilbert space $\ell {}^2 (V(G))$, we consider a von Neumann algebra. 
Let ${\cal B} (\ell {}^2 (V(G)))$ be the ${\bf C}^*$-algebra of bounded 
linear operators on $\ell {}^2 (V(G))$. 
A bounded linear operator $A$ of ${\cal B} (\ell {}^2 (V(G)))$ acts on 
$\ell {}^2 (V(G))$ by 
\[
A(f)(v)= \sum_{w \in V(G)} A(v,w)f(w), \   v \in V(G), \  f \in \ell {}^2 (V(G)) . 
\]
Then the von Neumann algebra ${\cal N}_0 (G, \Gamma )$ of bounded operators on 
$\ell {}^2 (V(G))$ commuting with the action of $\Gamma $ is defined as follows: 
\[
{\cal N}_0 (G, \Gamma )= \{ \lambda {}_0 (\gamma ) \mid \gamma \in \Gamma \} {}^{\prime} 
= \{ T \in {\cal B} (\ell {}^2 (V(G))) \mid \lambda {}_0 (\gamma )T=T \lambda {}_0 (\gamma ), 
\forall \gamma \in \Gamma \} . 
\]
The von Neumann algebra ${\cal N}_0 (G, \Gamma )$ inherits a trace 
by 
\[
{\rm Tr}_{\Gamma} (A)=\sum_{x \in {\cal F}_0 } \frac{1}{ \mid \Gamma {}_x \mid } A(x,x), 
\  A \in {\cal N}_0 (G, \Gamma ) . 
\]
Let the adjacency matrix ${\bf A} = {\bf A} (G)$ of $G$ be defined by 
\[
( {\bf A} f)(v)= \sum_{(v,w) \in R(G)} f(w), \  f \in \ell {}^2 (V(G)) . 
\]
By [14,15], we have 
\[
\mid \mid {\bf A} \mid \mid \leq d= \sup_{v \in V(G)} \deg {}_G v < \infty , 
\]
and so ${\bf A} \in {\cal N}_0 (G, \Gamma )$.

Similarly to $\ell {}^2 (V(G))$, we consider the Hilbert space 
$\ell {}^2 (R(G))$ of functions $f: R(G) \longrightarrow {\bf C} $ such that 
$\mid \mid \omega \mid \mid := \sum_{e \in R(G)} \mid \omega (e) \mid {}^2 
< \infty $. 
We define the left regular representation $\lambda {}_1 $ of $\Gamma $ on 
$\ell {}^2 (R(G))$ as follows: 
\[
(\lambda {}_1 ( \gamma ) \omega )(e)= \omega ( \gamma {}^{-1} e), 
\  \gamma \in \Gamma, \  \omega \in \ell {}^2 (R(G)), \  e \in R(G) . 
\]
Then the von Neumann algebra ${\cal N}_1 (G, \Gamma )= \{ \lambda {}_1 
(\gamma ) \mid \gamma \in \Gamma \} {}^{\prime} $ of bounded operators on 
$\ell {}^2 (R(G))$ commuting with the action of $\Gamma $, inherits a trace 
by 
\[
{\rm Tr}_{\Gamma} (A)=\sum_{e \in {\cal F}_1 } \frac{1}{ \mid \Gamma {}_e \mid } A(e,e), 
\  A \in {\cal N}_1 (G, \Gamma ) . 
\]

In an excellent paper [4], Fuglede and Kadison defined a positive-valued 
determinant for a von Neumann algebra with trivial center and finite trace $\tau $. 
For an invertible operator $A$ with polar decomposition $A=UH$, the 
Fuglede-Kadison determinant of $A$ is defined by 
\[
Det(A)= \exp \circ \tau \circ \log H, 
\]
where $\log H$ may be defined via functional calculus. 

Guido, Isola and Lapidus [6] extended the Fuglede-Kadison determinant 
to a determinant which is an analytic function. 
Let $( {\cal A} , \tau )$ be a von Neumann algebra with a finite trace 
$\tau $. 
Then, for $A \in {\cal A} $, let 
\[
\det {}_{\tau} (A)= \exp \circ \tau \circ \log A, 
\]
where 
\[
\log (A):= \frac{1}{2 \pi i} \int_{\Lambda } \log \lambda 
( \lambda -A )^{-1} d \lambda , 
\]
and $\Lambda $ is the boundary of a connected, simply connected region 
$\Omega $ containing the spectrum $\sigma (A)$ of $A$. 
Then the following lemma holds (see [7, Lemma 5.1]).

\section{A zeta function with respect to a generalization of a general coined quantum walk of an infinite periodic graph}

We define a zeta function with respect to a generalization of a general coined quantum walk of an infinite periodic graph. 

Let $G$ be a periodic graph with a countable discrete subgroup $\Gamma $ of $Aut \  G$. 
Moreover, let 
\[
{\bf I}_V =Id_{\ell {}^2 (V(G))} , {\bf I}_R =Id_{\ell {}^2 (R(G))} . 
\] 
Then, let ${\bf d}_i : \ell {}^2 (V(G)) \longrightarrow \ell {}^2 (R(G)) \ (i=1,2)$ such that 
\[
{\bf d}_i {\bf d}^*_i ={\bf I}_V . 
\] 
Furthermore, let 
\[
{\bf C}_i = a_i {\bf d}^*_i {\bf d}_i +b_i ( {\bf I}_R - {\bf d}^*_i {\bf d}_i ) \ (i=1,2) 
\] 
and ${\bf U} = {\bf C}_1 {\bf C}_2 $. 

Suppose that 
\[
{\bf C}^2_1 = {\bf I}_R . 
\]
Then we have 
\[
a_1= \pm 1 , \ b_1 = \mp 1. 
\]
Now, let 
\[
a_1 =1 , b_1 =-1. 
\]
A zeta function with respect to a general coined quantum walk of $G$ is defined 
as follows: 
\[
\zeta {} (G, \Gamma, u)= \det {}_{\Gamma } ( {\bf I}_R -u {\bf U} )^{-1} 
= \det {}_{\Gamma } ( {\bf I}_R -u ( a_1 {\bf d}^*_1 {\bf d}_1 + b_1 ( {\bf I}_R - {\bf d}^*_1 {\bf d}_1 ))
(a_2 {\bf d}^*_2 {\bf d}_2 +b_2 ( {\bf I}_R - {\bf d}^*_2 {\bf d}_2 )) )^{-1} . 
\]

Then we have the following result.

\begin{theorem}
Let $G$ be a periodic graph with a countable discrete subgroup $\Gamma $ 
of $Aut \  G$. 
Then 
\[
\det {}_{\Gamma } ({\bf I}_R -u {\bf U} )=(1+ b_2 u)^{{\rm Tr}_{\Gamma } ({\bf I}_R)-2{\rm Tr}_{\Gamma } ({\bf I}_V) } 
\det {}_{\Gamma } ((1+ a_2 u)(1-b_2 u) {\bf I}_V -2 c_2 u {\bf d}_2 {\bf d}^*_1 {\bf d}_1 {\bf d}^*_2 ) ,  
\]
where ${\rm Tr}_{\Gamma } ({\bf I}_R)= \sum_{e \in {\cal F}_1 } \frac{1}{ \mid \Gamma {}_e \mid } $ and 
$ {\rm Tr}_{\Gamma } ({\bf I}_V)= \sum_{v \in {\cal F}_0 } \frac{1}{ \mid \Gamma {}_v \mid } $(see [2]). 
\end{theorem}

{\bf Proof}.  The argument is an analogue of the method of Bass [1]. 

Let $G$ be a periodic graph with a countable discrete subgroup $\Gamma $ 
of $Aut \  G$. 

Now we consider the direct sum of the unitary representations 
$\lambda {}_0 $ and $\lambda {}_1 $: 
$\lambda (\gamma ):= \lambda {}_0 (\gamma ) \oplus \lambda {}_1 (\gamma )$  
$ \in {\cal B} (\ell {}^2 (V(G)) \oplus \ell {}^2 (R(G)))$. 
Then the von Neumann algebra $\lambda (\Gamma )^{\prime} 
:=\{ S \in {\cal B} (\ell {}^2 (V(G)) \oplus \ell {}^2 (R(G))) \mid 
S \lambda ( \gamma )= \lambda (\gamma )S, \gamma \in \Gamma \}$ 
consists of operators 
\[
S= 
\left[
\begin{array}{cc}
S_{00} & S_{01} \\
S_{10} & S_{11} 
\end{array}
\right]
,
\]
where $S_{ij} \lambda {}_j (\gamma )=\lambda {}_i (\gamma )S_{ij} , 
\gamma \in \Gamma , i,j=0,1$, so that $S_{ii} \in \Lambda {}_i 
\equiv  {\cal N}_i (G, \Gamma ), i=0,1$. 
Thus, $\lambda ( \Gamma )^{\prime} $ inherits a trace given by 
\[
{\rm Tr }_{\Gamma } 
\left[
\begin{array}{cc}
S_{00} & S_{01} \\
S_{10} & S_{11} 
\end{array}
\right]
:= {\rm Tr}_{\Gamma } (S_{00} )+ {\rm Tr}_{\Gamma } (S_{11} ) . 
\]

We introduce two operators as follows: 
\[
{\bf L} =
\left[
\begin{array}{cc}
(1- b^2_2 u^2 ) {\bf I}_V & -c_2 {\bf d}_2 -b_2 c_2 u {\bf d}_2 {\bf C}_1 \\
0 & {\bf I}_R 
\end{array}
\right]
,
{\bf M} = 
\left[
\begin{array}{cc}
{\bf I}_V & c_2 {\bf d}_2 +b_2 c_2 u {\bf d}_2 {\bf C}_1 \\
u {\bf C}_1 {\bf d}^*_2 & (1-b^2_2 u^2 ) {\bf I}_R 
\end{array}
\right]
,  
\]  
where $c_2=a_2-b_2$. 
Then we have 
\[
\begin{array}{rcl}
{\bf LM} & = & 
\left[
\begin{array}{cc}
(1- b^2_2 u^2 ) {\bf I}_V -c_2 u {\bf d}_2 {\bf C}_1 d^*_2 -b_2 c_2 u^2 {\bf d}_2 {\bf C}^2_1 {\bf d}^*_2 & 0 \\
u {\bf C}_1 {\bf d}^*_2 & (1- b^2_2 u^2 ) {\bf I}_R 
\end{array}
\right]
\\ 
\  &   &                \\ 
\  & = & 
\left[
\begin{array}{cc}
(1-a_2 b_2 u^2 ) {\bf I}_V -c_2 u {\bf d}_2 {\bf C}_1 {\bf d}^*_2 & 0 \\
u {\bf C}_1 {\bf d}^*_2 & (1- b^2_2 u^2 ) {\bf I}_R  
\end{array}
\right]
.
\end{array}
\]
Furthermore, we have 
\[
\begin{array}{rcl}
{\bf ML} & = & 
\left[
\begin{array}{cc}
(1- b^2_2 u^2 ) {\bf I}_V & 0 \\ 
u(1- b^2_2 u^2 ) {\bf C}_1 {\bf d}^*_2 & -c_2 u {\bf C}_1 {\bf d}^*_2 {\bf d}_2 - b_2 c_2 u^2 {\bf C}_1 {\bf d}^*_2 {\bf d}_2 {\bf C}_1 +(1- b^2_2 u^2 ) {\bf I}_R 
\end{array}
\right]
\\ 
\  &   &                \\ 
\  & = & 
\left[
\begin{array}{cc}
(1- b^2_2 u^2 ) {\bf I}_V & 0 \\ 
u(1- b^2_2 u^2 ) {\bf C}_1 {\bf d}^*_2 & ( {\bf I}_R -u(c_2 {\bf C}_1 {\bf d}^*_2 {\bf d}_2 +b_2 {\bf C}_1 ))( {\bf I}_R +u b_2 {\bf C}_1 ) 
\end{array}
\right]
.
\end{array}
\]
Here, note that ${\bf C}^2_1 = {\bf I}_R$.

For $\mid t \mid , \mid u \mid $ sufficiently small, we have 
\[
\sigma ( (1-a_2 b_2 u^2 ) {\bf I}_V -c_2 u {\bf d}_2 {\bf C}_1 {\bf d}^*_2 ), \sigma ((1- b^2_2 t^2 ) {\bf I}_V ), 
\sigma ((1- b^2_2 t^2 ) {\bf I}_R ), 
\]
\[
\sigma (( {\bf I}_R -u(c_2 {\bf C_1} {\bf d}^*_2 {\bf d}_2 +b_2 {\bf C}_1 ))( {\bf I}_R +u b_2 {\bf C}_1 ) )  
\in B_1 (1) = \{ z \in {\bf C} \mid \  \mid z-1 \mid <1 \} . 
\]
Similar to the proof of [7, Proposition 3.8], 
$\sigma ({\bf LM} )$ and $\sigma ({\bf ML})$ are contained in $B_1 (1)$. 
Thus, ${\bf L} $ and ${\bf M} $ are invertible, with bounded inverse, for 
$\mid t \mid , \mid u \mid $ sufficiently small. 

By Propositions 3.4, 3.6 and 3.8 in [7], we have 
\[
\begin{array}{rcl}
\det {}_{\Gamma } ( {\bf LM} ) & = & 
\det {}_{\Gamma } ((1- b^2_2 u^2 ) {\bf I}_V -c_2 u {\bf d}_2 {\bf C}_1 {\bf d}^*_2 -b_2 c_2 u^2 {\bf d}_2 {\bf C}^2_1 {\bf d}^*_2 ) 
\det {}_{\Gamma } ((1- b^2_2 u^2 ) {\bf I}_R) \\
\  &   &                \\ 
\  & = & 
(1- b^2_2 u^2 )^{{\rm Tr}_{\Gamma } ({\bf I}_R)} \det {}_{\Gamma } ((1- a_2 b_2 u^2 ) {\bf I}_V - c_2 u {\bf d}_2 {\bf C}_1 {\bf d}^*_2 ) 
\end{array}
\]
and 
\[
\begin{array}{rcl}
\det {}_{\Gamma } ( {\bf ML} ) & = & 
\det {}_{\Gamma } ((1- b^2_2 u^2 ) {\bf I}_V )
\det {}_{\Gamma } ( {\bf I}_R -u(c_2 {\bf C}_1 {\bf d}^*_2 {\bf d}_2 +b_2 {\bf C}_1 )) \det {}_{\Gamma } ( {\bf I}_R +u b_2 {\bf C}_1 ) \\
\  &   &                \\ 
\  & = & (1- b^2_2 u^2 )^{{\rm Tr}_{\Gamma } ({\bf I}_V)} 
\det {}_{\Gamma } ( {\bf I}_R -u(c_2 {\bf C}_1 {\bf d}^*_2 {\bf d}_2 +b_2 {\bf C}_1 )) \det {}_{\Gamma } ( {\bf I}_R +u b_2 {\bf C}_1 ) . 
\end{array}
\]

For $\mid t \mid , \mid u \mid $ sufficiently small, we have 
\[ 
{\bf ML}={\bf MLM}{\bf M}^{-1} ,
\]
and so, b[7, Proposition 3.7], 
\[ 
\det {}_{\Gamma } ( {\bf LM} ) = \det {}_{\Gamma } ( {\bf ML} ) . 
\]

Therefore, it follows that 
\[
(1- b^2_2 u^2 )^{{\rm Tr}_{\Gamma } ({\bf I}_R)} 
\det {}_{\Gamma } ( (1-a_2 b_1 u^2 ) {\bf I}_V -c_2 u {\bf d}_2 {\bf C}_1 {\bf d}^*_2 )
\]
\[
= (1- b^2_2 u^2 )^{ {\rm Tr}_{\Gamma } ({\bf I}_V)} 
\det {}_{\Gamma } ( {\bf I}_R -u {\bf C}_1 (c_2 {\bf d}^*_2 {\bf d}_2 +b_2 {\bf I}_R )) 
\det {}_{\Gamma } ( {\bf I}_R + b_2 u {\bf C}_1 ) , 
\]
and so 
\[
\begin{array}{rcl}
\  &  & \det {}_{\Gamma } ( {\bf I}_R -u {\bf C}_1 {\bf C}_2 )= \det {}_{\Gamma } 
( {\bf I}_R -u {\bf C}_1 (c_2 {\bf d}^*_2 {\bf d}_2 +b_2 {\bf I}_R )) \\
\  &   &                \\ 
\  & = & (1- b^2_2 u^2 )^{ {\rm Tr}_{\Gamma } ({\bf I}_R)
- {\rm Tr}_{\Gamma } ({\bf I}_V)} 
\det {}_{\Gamma } (  (1-a_2 b_2 u^2 ) {\bf I}_V -c_2 u {\bf d}_2 {\bf C}_1 {\bf d}^*_2 ) 
\det {}_{\Gamma } ( {\bf I}_R + b_2 u {\bf C}_1 )^{-1} . 
\end{array}
\]

But, we have 
\[ 
\begin{array}{rcl}
\  &  & \det {}_{\Gamma } ((1-a_2 b_2 u^2 ) {\bf I}_V -c_2 u {\bf d}_2 {\bf C}_1 {\bf d}^*_2 )  \\ 
\  &   &                \\ 
\  & = & \det {}_{\Gamma } ((1-a_2 b_2 u^2 ) {\bf I}_V -c_2 u {\bf d}_2 ( c_1 {\bf d}^*_1 {\bf d}_1 
+ b_1 {\bf I}_r ) {\bf d}^*_2 ) \\ 
\  &   &                \\ 
\  & = & \det {}_{\Gamma } ((1-b_1 c_2 u - a_2 b_2 u^2 ) {\bf I}_V 
-c_1 c_2 u {\bf d}_2 {\bf d}^*_1 {\bf d}_1 {\bf d}^*_2 ) ,  
\end{array}
\] 
where $c_1 = a_1 - b_1 =2$. 
Thus, we have 
\begin{equation} 
\begin{array}{rcl}
\det {}_{\Gamma } ( {\bf I}_R -u {\bf U} ) & = & (1- b^2_2 u^2 )^{{\rm Tr}_{\Gamma } ({\bf I}_R)- {\rm Tr}_{\Gamma } ({\bf I}_V)}  
\det {}_{\Gamma } ( {\bf I}_R + b_2 u {\bf C}_1 )^{-1} \\ 
\  &   &                \\ 
\  & \times & \det {}_{\Gamma } ((1-b_1 c_2 u - a_2 b_2 u^2 ) {\bf I}_V 
-c_1 c_2 u {\bf d}_2 {\bf d}^*_1 {\bf d}_1 {\bf d}^*_2 ) .  
\end{array}
\end{equation}

Next, we introduce two operators as follows: 
\[
{\bf L} =
\left[
\begin{array}{cc}
(1- b^2_1 u^2 ) {\bf I}_V & -c_1 {\bf d}_1 -b_1 c_1 u {\bf d}_1 \\
0 & {\bf I}_R 
\end{array}
\right]
,
{\bf M} = 
\left[
\begin{array}{cc}
{\bf I}_V & c_1 {\bf d}_1 +b_1 c_1 u {\bf d}_1 \\
u {\bf d}^*_1 & (1-b^2_1 u^2 ) {\bf I}_R 
\end{array}
\right]
,  
\]  
where $c_1=a_1-b_1$. 
Then we have 
\[
\begin{array}{rcl}
{\bf LM} & = & 
\left[
\begin{array}{cc}
(1- b^2_1 u^2 ) {\bf I}_V -c_1 u {\bf d}_1 d^*_1 -b_1 c_1 u^2 {\bf d}_1 {\bf d}^*_1 & 0 \\
u {\bf d}^*_1 & (1- b^2_1 u^2 ) {\bf I}_R 
\end{array}
\right]
\\ 
\  &   &                \\ 
\  & = & 
\left[
\begin{array}{cc}
(1-a_1 b_1 u^2 ) {\bf I}_V -c_1 u {\bf d}_1 {\bf d}^*_1 & 0 \\
u {\bf d}^*_1 & (1- b^2_1 u^2 ) {\bf I}_R  
\end{array}
\right]
.
\end{array}
\]
Furthermore, we have 
\[
\begin{array}{rcl}
{\bf ML} & = & 
\left[
\begin{array}{cc}
(1- b^2_1 u^2 ) {\bf I}_V & 0 \\ 
u(1- b^2_1 u^2 ) {\bf d}^*_1 & -c_1 u {\bf d}^*_1 {\bf d}_1 - b_1 c_1 u^2 {\bf d}^*_1 {\bf d}_1 +(1- b^2_1 u^2 ) {\bf I}_R 
\end{array}
\right]
\\ 
\  &   &                \\ 
\  & = & 
\left[
\begin{array}{cc}
(1- b^2_1 u^2 ) {\bf I}_V & 0 \\ 
u(1- b^2_1 u^2 ) {\bf d}^*_1 & ( {\bf I}_R -u(c_1 {\bf d}^*_1 {\bf d}_1 +b_1 {\bf I}_R ))( {\bf I}_R +u b_1 {\bf I}_R ) 
\end{array}
\right]
.
\end{array}
\]

For $\mid t \mid , \mid u \mid $ sufficiently small, we have 
\[
\sigma ( (1-a_1 b_1 u^2 ) {\bf I}_V -c_1 u {\bf d}_1 {\bf d}^*_1 ), \sigma ((1- b^2_1 u^2 ) {\bf I}_V ), 
\sigma ((1- b^2_1 u^2 ) {\bf I}_R ), 
\]
\[
\sigma (( {\bf I}_R -u(c_1 {\bf d}^*_1 {\bf d}_1 +b_1 {\bf I}_R ))( {\bf I}_R +u b_1 {\bf I}_R ) )  
\in B_1 (1) = \{ z \in {\bf C} \mid \  \mid z-1 \mid <1 \} . 
\]
Similar to the proof of [7, Proposition 3.8], 
$\sigma ({\bf LM} )$ and $\sigma ({\bf ML})$ are contained in $B_1 (1)$. 
Thus, ${\bf L} $ and ${\bf M} $ are invertible, with bounded inverse, for 
$\mid t \mid , \mid u \mid $ sufficiently small. 

By Propositions 3.4, 3.6 and 3.8 in [7], we have 
\[
\begin{array}{rcl}
\det {}_{\Gamma } ( {\bf LM} ) & = & 
\det {}_{\Gamma } ((1- b^2_1 u^2 ) {\bf I}_V -c_1 u {\bf d}_1 {\bf d}^*_1 -b_1 c_1 u^2 {\bf d}_1 {\bf d}^*_1 ) 
\det {}_{\Gamma } ((1- b^2_1 u^2 ) {\bf I}_R) \\
\  &   &                \\ 
\  & = & 
(1- b^2_1 u^2 )^{{\rm Tr}_{\Gamma } ({\bf I}_R)} 
\det {}_{\Gamma } ((1-a_1 b_1 u^2 ) {\bf I}_V -c_1 u {\bf d}_1 {\bf d}^*_1 ) 
\end{array}
\]
and 
\[
\begin{array}{rcl}
\det {}_{\Gamma } ( {\bf ML} ) & = & 
\det {}_{\Gamma } ((1- b^2_1 u^2 ) {\bf I}_V )
\det {}_{\Gamma } ( {\bf I}_R -u(c_1 {\bf d}^*_1 {\bf d}_1 +b_1 {\bf I}_R )) \det {}_{\Gamma } ( {\bf I}_R +u b_1 {\bf I}_R ) \\
\  &   &                \\ 
\  & = & (1- b^2_1 u^2 )^{{\rm Tr}_{\Gamma } ({\bf I}_V)} 
\det {}_{\Gamma } ( {\bf I}_R -u(c_1 {\bf d}^*_1 {\bf d}_1 +b_1 {\bf I}_R )) \det {}_{\Gamma } ( {\bf I}_R +u b_1 {\bf I}_R ) . 
\end{array}
\]

Let an orientation of $G$ be a choice of one oriented edge for each pair of 
edges in $R(G)$, 
which is called positively oriented. 
We denote by $E^+ G$ the set of positively oriented edges. 
Moreover, let $E^- G := \{ e^{-1} \mid e \in E^+ G \} $. 
An element of $E^- G$ is called a negatively oriented. 
Note that $R(G)=E^+ G \cup E^- G$. 

The operator ${\bf S}$ maps $\ell {}^2 (E^+ G)$ to $\ell {}^2 (E^- G)$. 
Then we obtain a representation $\rho $ of ${\cal B}( \ell {}^2 (R(G)))$ 
onto $Mat {}_2 {\cal B} ( \ell {}^2 (E^+ G))$, under 
\[
\rho ({\bf I}_R )=
\left[
\begin{array}{cc}
{\bf I} & 0 \\
0 & {\bf I} 
\end{array}
\right]
.
\]
By Propositions 3.6 and 3.8 in [7], 
\[
\det {}_{\Gamma } ({\bf I}_R +b_1 u {\bf I}_R )=  
\det {}_{\Gamma } 
\left[
\begin{array}{cc}
(1+ b_1 u) {\bf I} & 0 \\
0 & (1+ b_1 u) {\bf I} 
\end{array}
\right]
=(1+ b_1 u)^{ {\rm Tr}_{\Gamma} ( {\bf I}_R)} . 
\]

For $\mid t \mid , \mid u \mid $ sufficiently small, we have 
\[ 
{\bf ML}={\bf MLM}{\bf M}^{-1} ,
\]
and so, by [7, Proposition 3.7], 
\[ 
\det {}_{\Gamma } ( {\bf LM} ) = \det {}_{\Gamma } ( {\bf ML} ) . 
\]

Therefore, it follows that 
\[
(1- b^2_1 u^2 )^{{\rm Tr}_{\Gamma } ({\bf I}_R)} 
\det {}_{\Gamma } ( (1-a_1 b_1 u^2 ) {\bf I}_V -c_1 u {\bf d}_1 {\bf d}^*_1 )
\]
\[
= (1- b^2_1 u^2 )^{ {\rm Tr}_{\Gamma } ({\bf I}_V)} (1+ b_1 u)^{ {\rm Tr}_{\Gamma} ( {\bf I}_R)} 
\det {}_{\Gamma } ( {\bf I}_R -u (c_1 {\bf d}^*_1 {\bf d}_1 +b_1 {\bf I}_R )) , 
\]

But, we have 
\[ 
\begin{array}{rcl}
\  &   & \det {}_{\Gamma } ( (1-a_1 b_1 u^2 ) {\bf I}_V -c_1 u {\bf d}_1 {\bf d}^*_1 )= \det {}_{\Gamma } ( (1- c_1 u -a_1 b_1 u^2 ) {\bf I}_V ) \\
\  &   &                \\ 
\  & = & (1- c_1 u -a_1 b_1 u^2 )^{ {\rm Tr}_{\Gamma } ({\bf I}_V)} 
=(1- a_1 u)^{ {\rm Tr}_{\Gamma } ({\bf I}_V)} (1+ b_1 u)^{ {\rm Tr}_{\Gamma } ({\bf I}_V)} . 
\end{array}
\] 
Thus, we have 
\begin{equation} 
\begin{array}{rcl}
\  &  & \det {}_{\Gamma } ({\bf I}_R -u {\bf C}_1 ) \\ 
\  &   &                \\ 
\  & = & (1- b^2_1 u^2 )^{{\rm Tr}_{\Gamma } ({\bf I}_R)-{\rm Tr}_{\Gamma } ({\bf I}_V) } (1- a_1 u)^{ {\rm Tr}_{\Gamma } ({\bf I}_V)} 
(1+ b_1 u)^{ {\rm Tr}_{\Gamma } ({\bf I}_V)}  (1+ b_1 u)^{- {\rm Tr}_{\Gamma} ( {\bf I}_R)} \\ 
\  &   &                \\ 
\  & = & (1- b_1 u)^{{\rm Tr}_{\Gamma } ({\bf I}_R)-{\rm Tr}_{\Gamma } ({\bf I}_V) } (1- a_1 u)^{ {\rm Tr}_{\Gamma } ({\bf I}_V)} . 
\end{array}
\end{equation}
Substituting $-b_2 u$ into $u$ in (2), we have 
\[
\det {}_{\Gamma } ({\bf I}_R + b_2 u {\bf C}_1 )=(1+ b_1 b_2 u)^{{\rm Tr}_{\Gamma } ({\bf I}_R)-{\rm Tr}_{\Gamma } ({\bf I}_V) } 
(1+ a_1 b_2 u)^{ {\rm Tr}_{\Gamma } ({\bf I}_V)} . 
\] 
Hence, by (1), 
\[
\begin{array}{rcl}
\  &  & \det {}_{\Gamma } ({\bf I}_R -u {\bf U} ) \\ 
\  &   &                \\ 
\  & = & (1- b^2_2 u^2 )^{{\rm Tr}_{\Gamma } ({\bf I}_R)-{\rm Tr}_{\Gamma } ({\bf I}_V) } 
(1+ b_1 b_2 u)^{- {\rm Tr}_{\Gamma } ({\bf I}_R)+{\rm Tr}_{\Gamma } ({\bf I}_V) } 
(1+ a_1 b_2 u)^{- {\rm Tr}_{\Gamma } ({\bf I}_V)} \\ 
\  &   &                \\ 
\  & \times & \det {}_{\Gamma } ((1-b_1 c_2 u - a_2 b_2 u^2 ) {\bf I}_V 
-c_1 c_2 u {\bf d}_2 {\bf d}^*_1 {\bf d}_1 {\bf d}^*_2 ) \\ 
\  &   &                \\ 
\  & = & (1+ b_2 u)^{{\rm Tr}_{\Gamma } ({\bf I}_R)-2{\rm Tr}_{\Gamma } ({\bf I}_V) } 
\det {}_{\Gamma } ((1+ a_2 u)(1-b_2 u) {\bf I}_V -2 c_2 u {\bf d}_2 {\bf d}^*_1 {\bf d}_1 {\bf d}^*_2 ) . 
\end{array}
\] 
$\Box$

\end{document}